\documentclass[12pt]{article}
\usepackage{graphicx,color,microtype,breakcites}
\usepackage{amssymb,amsmath}

\newtheorem{Exa}{Example}
\newtheorem{thm}{Theorem}
\newtheorem{lem}{Lemma}
\newtheorem{cor}{Corollary}

\newtheorem{prop}{Proposition}
\newtheorem{Def}{Definition}
\newtheorem{Rem}{Remark}
\newtheorem{Exe}{Exercice}
\def\beXe{\begin{Exe}} \def\eeXe{\end{Exe}}
\def\eeD{\end{Def}} \def\beD{\begin{Def}}
\def\beXa{\begin{Exa}} \def\eeXa{\end{Exa}}
\def\beR{\begin{Rem}} \def\eeR{\end{Rem}}
\def\beL{\begin{Lem}} \def\eeL{\end{Lem}}
\newcommand{\finish}{\hfill$\Box$\vspace{0.2cm}}

\newcommand{\prf}{\noindent{\bf Proof:\ }}

\newcommand{\E}{{\rm I \!E}}
\newcommand{\p}{{\rm I \!P}}

\def\lev{L\'evy }      
 \def\T{T} \def\I{\infty} 
  \def\g{\gamma}    
  \def\th{\theta}
  
 \def\R{{\mathbb R}} \def\C{{\mathbb C}}
\def\BEN{\begin{enumerate}}  \def\BI{\begin{itemize}}
\def\EEN{\end{enumerate}}   \def\EI{\end{itemize}}
\def\be{\begin{equation}} \def\ee{\end{equation}}

\def\beq{\begin{eqnarray}}\def\eeq{\end{eqnarray}}
\def\bea{\begin{eqnarray*}}
\def\eea{\end{eqnarray*}}
\def\T{\widetilde}   \def\la{\label} 
\def\mD{\mathcal D} \def\fe{for example } 
\def\sn{spectrally negative }  \def\Y{Y}
\def\ol{\overline} \def\und{\underline}   \def\im{\item } \def\eqr{\eqref}  \def\wk{well known } \def\bc{\begin{cases}
  }   \def\dd{drawdown }
\def\ec{\end{cases}} \def\G{\Gamma} \def\Thr{Therefore, }
\long\def\symbolfootnote[#1]#2{
\begingroup
\def\thefootnote{\fnsymbol{footnote}}\footnote[#1]{#2}
\endgroup}
\def\fn{\symbolfootnote}

\providecommand{\keywords}[1]{\textbf{Keywords} #1}
\begin{document}




\title{Maximum Drawdown and Drawdown Duration of Spectrally Negative L\'{e}vy Processes Decomposed at Extremes  \thanks{This work is supported by T\"{u}bitak Project  with number 117F273}}

\author{Ceren Vardar-Acar \and Mine \c{C}a\u{g}lar \and Florin Avram}



\maketitle


\begin{abstract}  
	Path decomposition is performed to characterize the law of the pre/post-supremum, post-infimum and the intermediate processes of a spectrally negative L\'evy process taken up to an independent exponential time $T.$ As a result, mainly the distributions of the supremum of the post-infimum process and the maximum drawdown of the pre/post-supremum, post-infimum processes and the intermediate processes are obtained together with the law of drawdown durations.
\end{abstract}

\section{Introduction}
\label{intro}

Motivated by applications in statistics, mathematical finance
and risk theory, there has been increased interest recently in the study
of the running supremum and of the drawdown/regret/loss/process reflected at the supremum of a process $X$,
defined by
\[
\Y_t=\ol X_t-X_t, \quad \ol X_t:= \sup_{0 \leq s\leq t} X_s.
\]
For references to the numerous applications of drawdowns, see for example  \cite{MijoPist}, \cite{baur},  \cite{LLZU}, \cite{LLZ17}. Our  aim  is to find the distribution of the maximum loss/drawdown and the drawdown duration over various intervals with respect to the extremes of the process. In particular,
we generalize several results of \cite{salm} on Brownian motion.

The maximum loss/drawdown is defined by
\be M_t^ - : = \mathop {\sup }\limits_{0 \le u \le v \le t} (X_u^{} - X_v^{}) \ee
More generally, we will denote by
\be M_{a,b}^ - : = \mathop {\sup }\limits_{a \le u \le v \le b} (X_u^{} - X_v^{})  \ee
the corresponding quantities over an arbitrary interval.
We denote the first passage times above and below $x$ respectively by
\begin{equation*} \tau^{+}_x = \inf\{t\geq 0 : X_t>x\} ~~~~~~\tau^{-}_x = \inf\{t\geq 0 : X_t<x\} \;.\end{equation*}
Let the first passage time at level $d$ of the drawdown process be defined by
\be \label{al}
\alpha_d=\inf\{t\ge 0: Y_t > d\}
\ee
We define $\kappa_{\alpha_d}=\inf\{0\le s\le \alpha_d:\ol X_s =\ol X_{\alpha_{d}}\}$ as for spectrally negative Levy processes $\ol X$ is continuous. Then, the time
\be \label{dur}
\mathcal{T}^d = \alpha_{d}-\kappa_{\alpha_d}
\ee
defines the drawdown duration at the first time the drawdown exceeds level $d$.

We will develop below
formulas for maximum drawdown and drawdown duration under the assumption that $X$ is a
spectrally negative L\'evy processes  (with no positive jumps).
This class is popular   in numerous applications in risk theory, mathematical finance, queuing theory, etc, due to its combination of  flexibility and tractability via the usage of the so-called scale functions/$\gamma$-harmonic functions   $W^{(\gamma)}$ and  $Z^{(\gamma)}$ -- see \fe \cite{AKP,Kyp} and the recent cookbook \cite{AGV}.
In this paper, we add the latter cookbook more sophisticated results obtained by exploiting path decomposition into  the pre-supremum, post-supremum, post-infimum and  intermediate processes.


We perform several path decompositions through the extremes and characterize the distributions of the pre/post-supremum and post-infimum processes as Doob-$h$ transforms of the law of the L\'{e}vy process. These are given in Lemma 1 and Proposition 1, respectively. Moreover, conditioned on the event that the infimum occurs before the supremum, we also derive the $h$-transform for the intermediate process between the infimum and the supremum and post-supremum process in Proposition 2 in Section 3. These results lead to several useful distributions of maximum drawdown and drawdown duration given in Section 4. Explicitly, in Theorem 1, conditionally on the value of the supremum until an independent exponential time $T$ with parameter $\gamma$, denoted by $S_T=b$, $b>0$, the distributions of the maximum drawdown of the pre-supremum process, $M_{0,H_{S}}$, and the  post-supremum process, $M_{H_{S},T}$, are provided as
\[
\p\{M_{0,H_{S}}^{-}<d|S_T=b\} =e^{-b\left(\frac{ {W}^{ (\gamma)\prime}(d)}{ {W}^{ (\gamma)}(d)} -\Phi(\gamma)\right)}
\]
and
\begin{eqnarray*}
	\p\{M_{H_{S},T}^->d|S_T=b\}&=&\mathbb{E}[e^{ -\gamma \mathcal{T}^d_{H_{S},T}} | S_T=b]\\
	&=&  1+
	\frac{(Z^{(\gamma)}(d)-1)(W^{(\gamma)\prime}(d)-(W^{(\gamma)}(d))^2}
	{\Phi(\gamma) W^{(\gamma)}(d)}
\end{eqnarray*}
for $b >0$, $d>0$.
In addition, in Theorem 1, conditionally on the value of the infimum until  $T$ denoted by $I_T=b$, $a<0$,  the distribution of maximum drawdown of the post-infimum process, $M_{H_{S},T}$, is given by
\be \p\{M_{H_{I},T}^- > d |I_T=a\}=\p\{M_{H_{I},T}^-> d\}=1-\Phi(\gamma)	 \frac{W^{(\gamma)}(d)}{W^{(\gamma)^{\prime}}(d)}\ee for $d>a.$

In Theorem 2, we condition on the values of both extremes until $T$, including the infimum occurs at $a<0$, denoted as $I_T=a$, and the event that the infimum occurs before the supremum, denoted as $H_I<H_S$, and to be defined precisely below. We find the conditional distributions of the maximum drawdown of the intermediate process, denoted by $M_{H_I,H_{S}}$, and the post-supremum process as
\begin{eqnarray*}
	\lefteqn {\p\{M_{H_I,H_{S}}^-<d|H_I<H_S,I_T=a,S_T=b\}}\\
	& \displaystyle{ =
		\frac{  {W}^{(\gamma)}(b-a)}{{W}^{(\gamma)}(d)} e^{-(b-a-d)\left(\frac{{W}^{(\gamma)'}(d)}{{W}^{(\gamma)}(d)} -\Phi(\gamma){W}^{(\gamma)}(d)+\Phi(\gamma)\right)} }
\end{eqnarray*}
and
\begin{eqnarray*}\lefteqn {\p\{M_{H_{S},T}^-<d|H_I<H_S,I_T=a,S_T=b\}}\\
	&=&\mathbb{E}[e^{ -\gamma \mathcal{T}^d_{H_{S},T}} | H_I<H_S,I_T=a,S_T=b]\\
	&=&\frac{(Z^{(\gamma)}(d)-1)
		\frac{W^{ (\gamma)\prime}(d)}{W^{(\gamma)}(d)} - \gamma W^{(\gamma)}(d) }{  (Z^{(\gamma)}(b-a)-1)
		\frac{W^{ (\gamma)\prime}(b-a)}{W^{(\gamma)}(b-a)} - \gamma W^{(\gamma)}(b-a)   }
\end{eqnarray*}
for $0<d<b-a$ where the distrubutions of the drawdown duration for the post-supremum processes are found as corollaries.

Moreover, our results lead to several interesting observations related to  the $W,Z$ paradigm, which we expect to be useful elsewhere as well. Proposition 1 includes the law of the supremum of the post-infimum process, implying that
$$d \rightarrow \frac{\Phi(\gamma)({Z^{(\gamma)}}(d)-1)}{\gamma W^{(\gamma)}(d)}$$
for $d>0$, is a cumulative distribution function,  a fact that seems to have been overlooked before. This quantity could be important in applications, since it models the law
of a bull market (supremum) following a bear market (infimum).

Additionally, considering the process up to $\alpha_d$ rather than an independent exponential time $T$, in Proposition \ref{prop3} we identify the distribution of the part of the process after $\kappa_{\alpha_d}$ until $\alpha_d$ and find the Laplace transform of drawdown duration as
\[
\mathbb{E}[e^{-\gamma \mathcal{T}^{d}}]=\frac{W(d)}{W^{'}(d)}\frac{Z^{(\gamma)}(d)W^{(\gamma)'}(d)-\gamma W^{(\gamma)}(d)^2}{W^{(\gamma)}(d)} \; .
\]
This matches with the result given in \cite[Thm 3.1]{LLZ17}.

Independence of the decomposed processes  at the infimum or supremum for general Markov processes goes back the work of  \cite{mill}, and  may be  established directly for discrete time random walks by applying the strong Markov property  at weak ladder epochs (the successive times at which $X$ attains
an  infimum or supremum) \cite[Thm. 2.3]{bertoin1993splitting}. However, in the continuous time case, the proofs  require typically excursion theory  depending on whether $0$ is irregular or not for the reflected process; see for example \cite[Lem.VI.6.ii]{Ber}. \cite{chaum} specifies the laws of pre-infimum and post-infimum processes corresponding to a L\'{e}vy process conditioned to stay positive which starts at $x>0$, further properties of which are studied in \cite{chaumDon}. The latter process is characterized as an $h$-transform described with the local time of a L\'{e}vy process, for which 0 is assumed to be regular for $(-\infty, 0)$. As further related work, \cite{duquesne2003path} includes path decomposition for general Levy processes during an excursion of the drawdown process by splitting at the supremum of the excursion. The laws of the pre-supremum and post-supremum within an excursion are related to the process conditioned to stay positive and to that conditioned to stay negative, respectively. Similarly, the formulas which we derive in Lemma 1 are based on the observation of such conditioning, but for the path split at the supremum during a horizon of independent exponential time. In this paper, we follow \cite{chaum,chaumDon,duquesne2003path} and restrict ourselves to the unbounded variability case, when $0$ is regular for both $(-\infty,0)$ and $(0,\infty)$. The bounded variability case is left for future work.

The paper is organized as follows. The basics of spectrally negative L\'{e}vy processes are reviewed  in Section \ref{s:mod}. New path decomposition results for a spectrally negative L\'{e}vy process before an independent exponential time are given in Section \ref{s:pd}, and compared with some well-known probability laws. Our main results on the distribution of the maximum drawdown for the pre/post-supremum under two different sets of conditions, and the intermediate  processes  are presented in Section \ref{s:res}, together with results on the drawdown duration when applicable. Section \ref{s:prf} is devoted to the proof of the basic path decomposition result in Lemma 1 given in Section \ref{s:pd}.

\section{The \sn \lev model \la{s:mod}}

Let $(\Omega, {\cal H}, \p)$ be a probability space. We will assume  that under $\p$, $X=\{X_t\in \R, t\geq 0\}$ is a
L\'evy process, which  may therefore be characterized by its  L\'evy-Khinchine/Laplace exponent/symbol $\psi(\th)$   defined by
\beq \E_0 \left[e^{\th X_t }\right]=e^{t \; \psi(\th)}, \th \in \mD \subset \C, \la{LE}\eeq
where  $\mD$  includes at least the imaginary axis.

Furthermore, we assume that  $X$ is spectrally negative, which implies that its  Laplace exponent   may be written  as
\[
\psi (\lambda ) =  \mu \lambda  + \frac{{{\sigma ^2}}}{2}{\lambda ^2} + \int\limits_{( - \infty ,0)}^{} {({e^{\lambda x}} - 1 - \lambda x\;1_{\{x > -1\}})\, \Pi (dx)}
\]
where $\mu \in \mathbb{R}$, and $\sigma \geq 0$,
and where the L\'evy measure $\Pi$   is concentrated on $(-\infty, 0)$. 
Note that Brownian motion with drift $\mu$ can be recovered as a special case when $\Pi\equiv 0$. Using a Brownian motion $B$ and a Poisson random measure $N$, we can write
\[
dX_t=\mu\, dt+\sigma\, dB_t + \int_{(-1,0)} y\, \widetilde{N}(dy,dt) + \int_{(-\infty,-1]}y\, N(dy,dt)
\]
where $\widetilde{N}=N-\Pi(dy)\, dt$. We assume $\int_ {-\infty}^0 (1 \wedge y^2 ) \, \Pi(dy)< \infty$.
The infinitesimal generator  of $X$  is given by
\[
\mu f^\prime(x)+\frac{1}{2}\sigma^2  f^{\prime\prime}(x)+ \int_{-\infty}^0 [f(x+y)-f(x) -f^\prime(x)y\,1_{\{y>-1\}}]\, \Pi(dy)
\]
for  $f\in C_b^2$.

We will consider the process $X$ taken up to an exponential time $T$ independent of $X$, with parameter $\gamma.$
The running supremum and infimum until   time $T$ will be denoted by $S_T$ and $I_T$, respectively. Explicitly,
\be S_T:=\ol X_{T}=\sup \{ X_s: ~0  \leq s \leq T\}, I_T:=\und X_{T}:=\inf \{ X_s: ~0  \leq s \leq T\}. \ee
By \cite[Prop.2.2]{mill} the maximum and the minimum before a fixed time occur at a single point in time almost surely. The pre/post-supremum and  pre/post-infimum processes emerging thus  may be defined as follows.
\beD
Let
\be \la{SI}
H_S:={\inf}\{t < T: X_t=\ol X_t\} \quad \mbox{and} \quad H_I:={\inf}\{t < T:X_t =\und X_t\} \;.
\ee

The pre/post-$H_S$ and pre/post-$H_I$ processes are defined respectively by $$\bc X_u:0\leq u\leq H_{S}\\X_{H_S+u}
:0\leq u\leq T-H_S \\X_u:0\leq u\leq H_{I}\\X_{H_I+u}:0\leq \emph{}u\leq T-H_I\ec.$$

They give rise to four new probability measures, which are h-transforms of the original measure.
\eeD

The fluctuation identities that we rely on involve a class of functions known as scale functions. The $\gamma-$scale function of $X$ satisfies
\[\int\limits_0^\infty  {{e^{ - \lambda x}}} W^{(\gamma)}(x)dx = \frac{1}{{\psi (\lambda )}-\gamma} ~~~~\mbox{for}~~~~\lambda>\Phi(\gamma)  \]
where $\Phi$  denotes the right inverse of $\psi$ and the second scale function is given by
\begin{equation}  \label{Z}
{Z^{(\gamma)}}(x) = 1 + \gamma\int_0^x {{W^{(\gamma)}}(y)dy}
\end{equation}
see e.g. \cite[Thm. 8.1]{Kyp}.

\section{Path decomposition at  extrema \la{s:pd}}
The results of this paper will be obtained by the so called path decomposition, more precisely by splitting paths at  extrema, which provides a natural and efficient approach for further results in this field \cite{bertoin1993splitting,Ber,chaum,chaumDon,duquesne2003path}. For one example, think of splitting the process into a first piece preceding its global infimum, followed by  a second piece preceding the subsequent global supremum (the ''intermediate process"), and by a third piece following the global supremum. We will use these decomposition results for studying the maximum drawdown and the drawdown duration.

We state now a lemma  including a classic (part a))\cite{mill} and new results
(parts b),c)) for the  pre-$H_S$ process $\{X_u:0\le u\le H_{S}\}$ and post-$H_S$ process $\{X_{H_S+u}:0\le u\le T-H_S\}$.
\begin{lem} \la{lem1} Conditionally on $S_T=b$, it holds that:
	
	a) The pre-$H_S$  process and the post-$H_S$ process are independent.
	
	b)  the law of pre-$H_S$ process is an  $h$-transform of the law of the spectrally negative L\'{e}vy process killed at $T\wedge \tau_b^+$, with
	\[
	h(x)= e^{-\Phi(\gamma)(b-x)},  \quad x \leq b .
	\]
	
	More precisely, it is a spectrally negative L\'evy process killed when it hits $b$, and governed by an Esscher   exponential change of measure $\p^{\Phi(\gamma)}$ of the original
	measure, with
	\begin{equation}    \label{change}
	\left. \frac{ d \p^{\Phi(\gamma)}}{d\p} \right|_{{\cal F}_t}=e^{\Phi(\gamma)X_t -\g t}
	\end{equation}  and with Laplace exponent
	\[
	\T{\psi}(\lambda)=\psi_{\Phi(\gamma)}(\lambda):=\psi(\lambda + \Phi(\gamma))-\gamma,  \quad \lambda\ge -\Phi(\gamma).
	\]
	
	c) the law of the post-supremum process  is a $h$-transform of  the spectrally negative L\'{e}vy process killed at $T\wedge \tau_b^+$ with
	
	\[
	h(x)=1-e^{-\Phi(\gamma)(b-x)}, \quad  x\le b .
	\]
	
	Equivalently, the post-$H_S$ process is equal in distribution to $b+\G$, where the law of $\Gamma$ is the $h$-transform of the law of a spectrally negative L\'{e}vy  process conditioned to stay negative and killed at $T\wedge \tau_0^+$.
\end{lem}

The proof of Lemma \ref{lem1} is provided  in Section \ref{s:prf}.
Parts b)-c), which furnish an  explicit
description of the pre/post-supremum processes  and which generalize   \cite[Thm. 3.2]{salm}, do not seem to be available in the previous literature.

The following Proposition identifies the law of the post-infimum process as well as the distribution of the supremum in this part of the path.

\begin{prop} \label{pro}
	Conditionally on $I_T=a$, the pre-$H_I$ and post-$H_I$ processes are independent, and the law of the post-$H_I$ process is the $h$-transform of the law of the spectrally negative L\'{e}vy process killed at $T\wedge \tau_a^-$ with
	\[
	h(x)=1- Z^{(\gamma)}(x-a)+\frac{\gamma}{\Phi(\gamma)}W^{(\gamma)}(x-a) \; .
	\]
	Moreover, the distribution of the supremum of the post-infimum process, $S_{H_{I},T}$, when $I_T=a$ is given by
	\be \la{IT}\p(S_{H_{I},T}\le b|I_T=a)=\frac{\Phi(\gamma)({Z^{(\gamma)}}(b-a)-1)}{\gamma W^{(\gamma)}(b-a)}\ee
	for $b>a$.
\end{prop}
\prf
The post-$H_I$ process is a Markov process with stationary transitions, which depend on the value of the infimum, $I_T$ \cite{mill}.
The post-$H_I$ process evolves as a spectrally negative L\'evy process which is conditioned to stay above $I_T$ and killed at the minimum of the exponential time $T$ and $\tau_a^-$. We can identify its transition function $P_t$ as
\begin{equation}
P_t(x,dy) = \p_x\{X_t\in dy, t<T\wedge \tau_a^- \,|\, T< \tau_a^-\}  \label{post}
\end{equation}
which is equal to
\begin{eqnarray*}
	& &  \displaystyle{\frac{  \p_x\{X_t\in dy, t<T \wedge \tau_a^-, T< \tau_a^-\} }{\p _x\{T< \tau_a^-\} }} \nonumber \\
	& & =\displaystyle{ \frac{  \p_x\{ T< \tau_a^- | X_t = y, t<T \wedge \tau_a^- \}\p _x  \left\{   X_t\in dy, t<T \wedge \tau_a^- \right\}}{\p _x\{T< \tau_a^-\} } } \nonumber \\
	& & = \displaystyle{\frac{  \p_y\{ T\circ \theta_t< \tau_a^-\circ \theta_t \}\p _x  \left\{   X_t\in dy, t<T \wedge \tau_a^- \right\}}{\p _x\{T< \tau_a^-\} } }\nonumber \\
	& & = \displaystyle{ \frac{  \p_{y} \{ T< \tau_a^-\} }{\p _x\{T< \tau_a^-\}} \p _x  \left\{   X_t\in dy, t<T \wedge \tau_a^- \right\} } \nonumber
\end{eqnarray*}
by Markov property and the fact that $T$ is memoryless, where $\theta$ is the shift operator \cite[Sec.7.3]{cinl}. Let ${h}(x) = \p _x\{T< \tau_0^-\}$, which can be evaluated as
\begin{eqnarray}
{h}(x) &=& \E_x[\p_x\{T<\tau_a^-|\, \tau_a^-\}] = \E_x[1-e^{-\gamma \tau_a^-}1_{\{\tau_a^-<\infty\}}] \nonumber \\
&=& 1-\E_x[e^{-\gamma \tau_a^-}1_{\{\tau_a^-<\infty\}}] \;=
\;  1-\E_{x-a}[e^{-\gamma \tau_0^-}1_{\{\tau_0^-<\infty\}}]  \nonumber \\
&=& 1- Z^{(\gamma)}(x-a)+\frac{\gamma}{\Phi(\gamma)}W^{(\gamma)}(x-a)\nonumber  \; .
\end{eqnarray}
Note that the transition density of the post-infimum process can be written as
\[
P_t(x,dy)=\frac{{h}(y)}{{h}(x)}\p_{x}\left\{ X_t\in dy, t<T\wedge\tau_a^- \right\}
\]
for $x,y>a$. Then, it follows that
\begin{eqnarray*} \p(S_{H_{I},T}>b|I_T=a)&=& \lim_{x \to a} \frac{1}{{h}(x)}\, \E_{x}  \left[ {h}(X_{\tau_{b}^+});\tau_{b}^+<\tau_a^-, \tau_{b}^+<T \right]\\
	&=&\lim_{x \to a} \frac{{h}(b)}{{h}(x)}\frac{W^{(\gamma)}(x-a)}{W^{(\gamma)}(b-a)}\\
	&=&\lim_{x \to a}\frac{\frac{\gamma W^{(\gamma)}(b-a)}{\Phi(\gamma)}-\gamma \int_0^{b-a} W^{(\gamma)}(z)\, dz }{\frac{\gamma W^{(\gamma)}(x-a)}{\Phi(\gamma)}-\gamma \int_0^{x-a} W^{(\gamma)}(z)\, dz }\: \frac{W^{(\gamma)}(x-a)}{W^{(\gamma)}(b-a)}\\
	&=&1-\frac{\Phi(\gamma)}{W^{(\gamma)}(b-a)}\int_0^{b-a} W^{(\gamma)}(z)\, dz
\end{eqnarray*}
which is simplified by \eqref{Z} and we have used the fact that $W^{(\gamma)}(0)=0$ because $X$ is assumed to be of unbounded variation.
\finish

As a special case, the  function ${h}(x)$  corresponding to post-$H_I$ process in Proposition \ref{prop3} is equal to $1-e^{-x\sqrt{2\gamma}}$ for standard Brownian motion due to symmetry and we recover the result in \cite[Thm.3.2.1]{salm}. Moreover, the above theorem is a generalization of \cite[Eq.3.23]{salm}.

\beR  Note that the law of the post-$H_I$ process is that of $a+\Gamma$, where $\Gamma$ has  the  law $\p$ of the original  process conditioned to stay positive  and killed at time $T\wedge \tau_0^-$. The law of  \lev processes conditioned to stay positive, denoted by $\p^\uparrow$, was constructed in \cite{bertoin1993splitting} having the $h$-function equal to the scale function $W$ \cite[pg.200]{Ber}. Observation of the processes over an exponential horizon of rate $\gamma$ changes the $h$-function as given in Proposition \ref{prop3}, and not simply to $W^{(\gamma)}$. Similarly,  the law of the post-$H_S$ process is related to the law of the original  process conditioned to stay negative  and killed at time $T\wedge \tau_0^+$ as given in Lemma \ref{lem1} with an $h$-function different from that corresponding to the law denoted by $\p^\downarrow$ (see e.g. \cite{duquesne2003path}).
\eeR

We now identify the conditional law of the decomposed paths between the infimum and the supremum, and post-supremum in the following proposition. We call the process $\{X_{H_I+u}: 0\le u\le H_S - H_I\}$ as the \emph{intermediate process}.

\begin{prop} \label{prop2} Conditionally on $H_I<H_S$, $I_T=a$, $S_T=b$, the intermediate process  and the post-$H_S$ process are independent, and have the following distributions:
	
	i. the law of the intermediate process is that of $a+\Gamma$, where the law of $\Gamma$ is the $h$-transform of the  spectrally negative L\'{e}vy process killed at $\tau_{b-a}^+$, with
	\[
	h(x) = e^{-\Phi(\gamma)(b-x)}\, \frac{W^{ (\gamma)}(x)}{W^{(\gamma)}(b)}
	\]
	
	ii. the law of the post-$H_S$ process is the $h$-transform of the law of the spectrally negative L\'{e}vy process killed at $T \wedge {\tau}^+_{b} \wedge {\tau}_a^-$ with
	\[
	h(x)=1-Z^{(\gamma)}(x-a)
	+ (Z^{(\gamma)}(b-a)-1)\frac{W^{ (\gamma)}(x-a)}{W^{(\gamma)}(b-a)}
	\]
\end{prop}
\prf i.  By Lemma \ref{lem1}, the pre-$H_S$ process is a L\'{e}vy process with Laplace exponent $\T{\psi}$. Let $\p^{\Phi(\gamma)}$ denote the probability corresponding to this change of measure.  Now, first note that the intermediate process being before the supremum is clearly governed by  $\p^{\Phi(\gamma)}$.  Next, the infimum of the presupremum process is also equal to $a$ under the given conditions.
Therefore, the intermediate process starts at $a$ and evolves independently from the past before the infimum as a process conditioned to stay above $a$ and killed at the first passage time above $b$ before the exponential time $T$. This implies that we can write the probability transition function of the intermediate process as
\begin{equation}   \label{11}
P_t(x,dy) =\p_x^{\Phi(\gamma)}\{X_t\in dy, t<\tau_b^+\,|\, \tau_{b}^+<\tau_{a}^-\wedge T\}
\end{equation}
In view of Lemma 1, in particular, \eqref{eqSec6} in its proof, we get
\begin{equation} \label{above}
P_t(x,dy) = \frac{e^{-\Phi(\gamma)(b-y)}}{e^{-\Phi(\gamma)(b-x)}} \p_x \{X_t\in dy, t<\tau_b^+\wedge T \,|\, \tau_{b}^+<\tau_{a}^-\wedge T\}
\end{equation}
We evaluate  the conditional distribution on the right hand side of \eqref{above} as
\begin{eqnarray*}
	&& \hspace{-1cm} \displaystyle{\frac{  \p_x\{X_t\in dy, t<\tau_{b}^+ , \tau_{b}^+<\tau_{a}^-\wedge T\} }{\p _x\{\tau_{b}^+<\tau_{a}^-\wedge T\} }}   \\
	& =& \displaystyle{ \frac{  \p_x\{ \tau_{b}^+<\tau_{a}^-\wedge T | X_t = y, t<\tau_{b}^+ \}\p _x  \left\{   X_t\in dy, t<\tau_{b}^+ \right\}}{\p _x\{\tau_{b}^+<\tau_{a}^-\wedge T\} } } \\
	& =& \displaystyle{\frac{  \p _y\{ \tau_{b}^+\circ \theta_t< \tau_a^-\circ \theta_t \wedge T \circ \theta_t \}\p _x  \left\{   X_t\in dy, t<\tau_{b}^+  \right\}}{\p _x\{\tau_{b}^+<\tau_{a}^-\wedge T\}} }  \\
	& =& \displaystyle{ \frac{  \p_{y} \{ \tau_{b}^+<\tau_{a}^-\wedge T\} }{\p _x\{\tau_{b}^+<\tau_{a}^-\wedge T\}} \:\p _x  \left\{   X_t\in dy, t< \tau_{b}^+  \right\} }
\end{eqnarray*}
by Markov property. Plugging this expression  in \eqref{above}, we have
\[
P_t(x,dy) = \frac{e^{-\Phi(\gamma)(b-y)}   \p_{y} \{ \tau_{b}^+<\tau_{a}^-\wedge T\} }
{e^{-\Phi(\gamma)(b-x)}\p _x\{\tau_{b}^+<\tau_{a}^-\wedge T\} }  \:\p _x \{   X_t\in dy, t< \tau_{b}^+  \} \; .
\]
As a result, the distribution of the intermediate process is Doob $h$-transform of the law of the original L\'{e}vy process killed at $\tau_b^+$ with
\begin{eqnarray}
h(x):&=&e^{-\Phi(\gamma)(b-x)}\p _x\{\tau_{b}^+<\tau_{a}^-\wedge T\} \label{**}\\
&=& e^{-\Phi(\gamma)(b-x)} \E _{x-a}[e^{-\gamma\tau_{b-a}^+};\tau_{b-a}^+ <\tau^-_{0}] \nonumber \\
&= & e^{-\Phi(\gamma)(b-x)}\, \frac{W^{ (\gamma)}(x-a)}{W^{(\gamma)}(b-a)}  \nonumber
\end{eqnarray}
by  the two sided exit formula   \cite[pg.24]{KKR}.


ii. Under the condition $I_T=a$, $S_T=b$, $H_I<H_S$, for $x<b$, $a<y<b$, the transition semigroups $P_t(x,dy)$ of the post-$H_S$ process can be found in terms of the transition probability density of the original process with an $h$-transform because starting at $x$, the process is conditioned to stay above  $a$ and below $b$, before an exponential time $T$.
As a result, the transition density of the post-$H_S$ process under the given conditions is the $h$-transform of the transition density of the process until $T \wedge {\tau}^+_{b} \wedge {\tau}_a^-$ with
\begin{eqnarray*}\label{h} h(x) &=&  \p_{x} \{ T<{\tau}^+_{b} \wedge{\tau}_a^-\}  \\
	&=& \p _{x}\{T< \tau^-_{a}, \tau^-_{a}<\tau_{b}^+\}+\p _{x}\{T<\tau_{b}^+, \tau_{b}^+ < \tau^-_{b}\} \\
	&=& \p_{x}\{\tau^-_{a}<\tau_b^+\}-\E _{x}[e^{-\gamma\tau^-_{a}};{ \tau^-_{a}<\tau_b^+}] \\
	&&\;+ \p _{x}\{\tau_b^+ < \tau^-_{a}\}-\E _{x}[e^{-\gamma\tau_b^+};\tau_b^+ <\tau^-_{a}]\\
	&=&  1-      \E _{x-a}[e^{-\gamma\tau^-_{0}};{ \tau^-_{0}<\tau_{b-a}^+}] -\E _{x-a}[e^{-\gamma\tau_{b-a}^+};\tau_{b-a}^+ <\tau^-_{0}]                         \\
	&=& 1-Z^{(\gamma)}(x-a)
	+ (Z^{(\gamma)}(b-a)-1)\frac{W^{ (\gamma)}(x-a)}{W^{(\gamma)}(b-a)}
\end{eqnarray*}
where   the two sided exit formulas given in \cite[pgs.17,24]{KKR} are used. \finish

\beR
Alternatively, the very same distribution   has been obtained by \cite{arxive} through enlargement of filtrations at the last exit from the infimum.
Although we consider the first exit from the infimum rather than the last exit, the two proofs are equivalent because the maximum and the minimum occur at a single point in time almost surely \cite[Prop.2.2]{mill}. The process from the last exit from the infimum, namely $a$,  until the first passage above $b$ turns out to be a diffusion with an infinitesimal generator revealing the very same $h$-transform to the law of the original L\'{e}vy process observed until $\tau_b^+$. Explicitly, the generator given in \cite[pg.17]{arxive} is
\begin{eqnarray*}
	{\mathcal L}F(x) & =&  \left[\mu +\sigma \frac{ {W}^{ (\gamma) \prime}(x)-\Phi(\gamma) {W}^{ (\gamma)}(x)}
	{ {W}^{(\gamma)}(x)}\right] F'(x)+\frac{\sigma^2}{2} F''(x) \\
	& & + \int_{(-1,0)}y \frac{e^{-\Phi(\gamma)y} {W}^{ (\gamma)}(x+y) -  {W}^{ (\gamma)}(x)} { {W}^{(\gamma)}(x) }      \,\Pi(dy) \\
	&& + \int_{-\infty}^0 [F(x+y)-F(x)-F'(x)\, y\, 1_{\{y>-1\}} ]     \\
	&& \hspace{1cm} .\,\left[1+ \frac{e^{-\Phi(\gamma)y} {W}^{ (\gamma)}(x+y) -  {W}^{ (\gamma)}(x)} { {W}^{(\gamma)}(x) } \right] \,\Pi(dy)
\end{eqnarray*}
which reveals that
\[
\frac{h'(x)}{h(x)} = \frac{ {W}^{ (\gamma) \prime}(x)-\Phi(\gamma) {W}^{ (\gamma)}(x)}
{ {W}^{(\gamma)}(x)}
\]
consistent with \eqref{**}, free from the start and end points, $a$ and $b$ respectively.
\eeR


\section{Maximum Drawdown and Drawdown Duration \la{s:res}}
We provide the distribution of the maximum drawdown for the pre-supremum and post-supremum processes using the results of the Section 3.
\begin{thm} \label{thm3} We have,
	\BEN \im conditionally on $S_T=b$,
	\[
	\p\{M_{0,H_{S}}^{-}<d|S_T=b\} =e^{-b\left(\frac{ {W}^{ (\gamma)\prime}(d)}{ {W}^{ (\gamma)}(d)} -\Phi(\gamma)\right)}
	\]
	\im conditionally on $S_T=b$,
	\be \la{MS}
	\p\{M_{H_{S},T}^->d|S_T=b\}
	= 1+
	\frac{(Z^{(\gamma)}(d)-1)W^{(\gamma)\prime}(d)-\gamma(W^{(\gamma)}(d))^2}
	{\Phi(\gamma) W^{(\gamma)}(d)}
	\ee for all $b >0, d>0.$
	\Thr this holds unconditionally as well:
	\bea
	\p\{M_{H_{S},T}^->d\}
	=  1+
	\frac{(Z^{(\gamma)}(d)-1)W^{(\gamma)\prime}(d)-\gamma(W^{(\gamma)}(d))^2}
	{\Phi(\gamma) W^{(\gamma)}(d)},
	\eea
	for all $d>0.$
	\im conditionally on $I_T=a$, which also holds unconditionally,
	\be \la{IMDD}\p\{M_{H_{I},T}^- > d |I_T=a\}=\p\{M_{H_{I},T}^-> d\}=1-\Phi(\gamma)	 \frac{W^{(\gamma)}(d)}{W^{(\gamma)^{\prime}}(d)}\ee for $d>a.$
	\EEN
\end{thm}

\prf i. Let $\T{W}(\cdot)= {W}_{\Phi(\gamma)}^{ (0)}(\cdot)$ denote the $0$-scale function corresponding to the Laplace exponent $\T{\psi}(\lambda)$ of the   Esscher transform $\p^{\Phi(\gamma)}$.
By Lemma \ref{lem1}b), we have
\begin{eqnarray*} \p\{M_{0,H_{S}}^{-}<d|S_T=b\} = \p^{\Phi(\gamma)}
	\{{M} _{{\tau}_b^+}^{-}<d, {\tau}_b^+<T\}  =
	e^{-b \frac{\T{W}^{ (\gamma)\prime}(d)}
		{\T{W}^{(\gamma)}(d)} }
\end{eqnarray*}
by \cite[Eq.13]{avram2019w} in view of
\[
\{ {M} _{\tau_b^+}^{-}<d\}= \{\tau_{b}^+<\alpha_{d}\}.
\]
Finally, the \wk Esscher transform identity ${W}_{\Phi(\gamma)}^{ (0)}(x)= e^{-\Phi(\gamma)x} {W}^{(\gamma)}(x)$ \cite[pg.236]{Kyp} yields the result.

ii. Now, for the maximum drawdown of the post-supremum process to be greater than $d>0$, we must have the passage time of the process below $b-d$ to occur before $T$.  Considering the condition $\{T<\tau_b^+\}$  as implied by $S_T=b$ and confirmed with the $h$-transform in Lemma \ref{lem1}c), and using the drawdown time $\alpha_d$, we get
\begin{eqnarray*} \p\{M_{H_{S},T}^- \!\!\!\!  \!\!\!\! \!\!\!\! &&> d |S_T=b\}
	\\
	&=& \lim_{x \to b-}\p_x\{\alpha_d<T |T<\tau_b^+\} = \lim_{x \to b-}\frac{\p_x\{\alpha_d<T ,T<\tau_b^+\}}{\p_x\{T<\tau_b^+\}} \\
	&=& \lim_{x \to b-}\frac{\p_x\{\alpha_d<T\}- \p_x\{\alpha_d<T,\tau_b^+<T\}}{\p_x\{T<\tau_b^+\}} \\
	&=& \lim_{x \to b-}\frac{\p_x\{\alpha_d<T\}-\p_x\{\alpha_d<T ,\tau_b^+<\alpha_d\}-\p_x\{\tau_b^+<T ,\alpha_d<\tau_b^+\}}{\p_x\{T<\tau_b^+\}}\\
	&=& \lim_{x \to b-}\frac{ \p_x\{\alpha_d<T,\alpha_d<\tau_b^+\}-\p_x\{\tau_b^+<T\}+ \p_x\{\tau_b^+<T,\tau_b^+<\alpha_d\}}{\p_x\{T<\tau_b^+\}}\\
	&=&\lim_{x \to b-}\frac{\left(1-e^{-(b-x) \frac{W^{(\gamma)^{\prime}}(d)}{W^{(\gamma)}(d)}}\right)\left(Z^{(\gamma)}(d)-W^{(\gamma)}(d)\frac{Z^{(\gamma)^{\prime}}(d)}{W^{(\gamma)^{\prime}}(d)}\right)-e^{-\Phi(\gamma)(b-x)}+e^{-(b-x) \frac{W^{(\gamma)^{\prime}}(d)}
			{W^{(\gamma)}(d)}}}{1-e^{-\Phi(\gamma)(b-x)}}\\
	&=&\frac{\left(Z^{(\gamma)}(d)-W^{(\gamma)}(d)
		\frac{Z^{(\gamma)^{\prime}}(d)}{W^{(\gamma)^{\prime}}(d)}\right)
		\frac{W^{(\gamma)^{\prime}}(d)}{W^{(\gamma)}(d)}-
		\frac{W^{(\gamma)^{\prime}}(d)}{W^{(\gamma)}(d)}+
		\Phi(\gamma)}
	{+\Phi(\gamma)}\\
	&=&  1+
	\frac{(Z^{(\gamma)}(d)-1)W^{(\gamma)\prime}(d)-\gamma(W^{(\gamma)}(d))^2}
	{\Phi(\gamma) W^{(\gamma)}(d)}
\end{eqnarray*}
using  equations (11), (13), and (15) in \cite{avram2019w} together with \eqref{ilkh}  above, and applying L'Hospital's rule once (the result was simplified using $Z^{(\gamma)\prime}= \gamma W^{(\gamma)}$).

iii. When $I_T=a$ is given by Proposition \ref{pro} we have,
\begin{eqnarray*}&& \p\{M_{H_{I},T}^- > d |I_T=a\}
	\\
	&=& \lim_{x \to a+}\p_x\{\alpha_d<T |T<\tau_a^-\} = \lim_{x \to a+}\frac{\p_x\{\alpha_d<T ,T<\tau_a^-\}}{\p_x\{T<\tau_a^-\}} \\
	&=& \lim_{x \to a+}\frac{\p_x\{\alpha_d<T\}- \p_x\{\alpha_d<T,\tau_a^-<T\}}{\p_x\{T<\tau_a^-\}} \\
	&=& \lim_{x \to a+}\frac{\p_x\{\alpha_d<T\}-\p_x\{\alpha_d<T ,\tau_a^-<\alpha_d\}-\p_x\{\tau_a^-<T ,\alpha_d<\tau_a^-\}}{\p_x\{T<\tau_a^-\}}\\
	&=& \lim_{x \to a+}\frac{\p_x\{\alpha_d<T,\alpha_d<\tau_a^-\}-\p_x\{\tau_a^-<T\}+\p_x\{\tau_a^-<T ,\tau_a^-<\alpha_d\}}{\p_x\{T<\tau_a^-\}}\\
	&=& \lim_{x \to a+}\frac{\left(Z^{(\gamma)}(d)-W^{(\gamma)}(d)
		\frac{Z^{(\gamma)^{\prime}}(d)}{W^{(\gamma)^{\prime}}(d)}\right)
		\frac{W^{(\gamma)}(x-a)}{W^{(\gamma)}(d)}+\frac{\gamma}{\Phi(\gamma)}W^{(\gamma)}(x-a)-\frac{W^{(\gamma)}(x-a)}{W^{(\gamma)}(d)}Z^{(\gamma)}(d)}{1-Z^{(\gamma)}(x-a)+\frac{\gamma}{\Phi(\gamma)}W^{(\gamma)}(x-a)}\\
	&=&\lim_{x \to a+}\frac{\frac{\gamma}{\Phi(\gamma)}W^{(\gamma)}(x-a)-\frac{\gamma W^{(\gamma)}(d)}{W^{(\gamma)^{\prime}}(d)}W^{(\gamma)}(x-a)}{1-Z^{(\gamma)}(x-a)+\frac{\gamma}{\Phi(\gamma)}W^{(\gamma)}(x-a)}\\
	&=&\frac{(\frac{\gamma}{\Phi(\gamma)}-\frac{\gamma W^{(\gamma)}(d)}{W^{(\gamma)^{\prime}}(d)})W^{(\gamma)^{\prime}}(0^+)}{\frac{\gamma}{\Phi(\gamma)}W^{(\gamma)^{\prime}}(0^+)}=1-\Phi(\gamma)	 \frac{W^{(\gamma)}(d)}{W^{(\gamma)^{\prime}}(d)}
\end{eqnarray*} using  equations (8) and \cite[Theorem 1, last row, middle column]{avram2019w} and applying L'Hospital's rule once (the result was simplified using $Z^{(\gamma)\prime}= \gamma W^{(\gamma)}$).
\finish

\beR The fact that \eqr{MS} and \eqr{IMDD} are independent of $b$ and $a,$ respectively is surprising at first, but in line with the \wk fact that the law of the \dd for a \lev process
is independent  of its starting point.\eeR

\beR \la{wow} It is easy to check that the function $d \rightarrow e^{-b\left(\frac{ {W}^{ (\gamma)\prime}(d)}{ {W}^{ (\gamma)}(d)} -\Phi(\gamma)\right)}$ is indeed a cdf since the  rate of excursions larger than $d$, given by
$\frac{ {W}^{ (\gamma)\prime}(d)}{ {W}^{ (\gamma)}(d)}$, decreases from $\I$ to $\Phi(\gamma)$. Similarly $d \rightarrow \Phi(\gamma)	\frac{W^{(\gamma)}(d)}{W^{(\gamma)^{\prime}}(d)}$ is a cdf.

The fact that the function $d \rightarrow 1+
\frac{(Z^{(\gamma)}(d)-1)(W^{(\gamma)\prime}(d)-(W^{(\gamma)}(d))^2}
{\Phi(\gamma) W^{(\gamma)}(d)}$ is a cdf is a new result and may be verified numerically.

\eeR

Now, the Laplace transform of the drawdown duration defined in  \eqref{dur} for a given $d>0$ becomes readily available for post-$H_S$ process from Theorem \ref{thm3}. To see this, note that
\[
\{\alpha_{d}(H_{S},T)  <T\}=\{ M_{H_{S},T}^->d \}
\]
conditionally on $S_T=b$, where
\[
\alpha_{d}(H_{S},T)=\inf \{ H_S<t<T: Y_t>d\}
\]
is the first passage of the drawdown of the post-$H_S$ process above $d$. Since the post-supremum process stays below $b$, the drawdown duration  during $(H_S,T)$, denoted by $\mathcal{T}^d_{H_{S},T} $ satisfies
\[
\mathcal{T}^d_{H_{S},T}  = \alpha_{d}(H_{S},T) \; .
\]
Then, Theorem \ref{thm3} provides the following formula.
\begin{cor} \label{cor1} We have
	\[
	\mathbb{E}[e^{ -\gamma \mathcal{T}^d_{H_{S},T}} | S_T=b] = \mathbb{E}[e^{ -\gamma \mathcal{T}^d_{H_{S},T}}]= 1-Z^{(\gamma)}(d)\frac{W^{(\gamma)^{\prime}}(d)}{\Phi(\gamma)W^{(\gamma)}(d)}
	-\frac{\gamma W^{(\gamma)}(d)}{\Phi(\gamma) }\; .
	\]
\end{cor}

Under the condition that the infimum of the L\'evy process occurs before its supremum with given levels, we  identify the distribution of the maximum drawdown of the intermediate process, as well as the post-supremum process.
The following theorem is for the distribution of the maximum drawdown for the intermediate and post-supremum processes.

\begin{thm}\label{thm2} Conditionally on $H_I<H_S$, $I_T=a$, $S_T=b$, it follows that \\
	
	i. the distribution of the maximum loss of the intermediate process, $M_{H_{I},H_S}^-$, when $H_I<H_S$, $I_T=a$, $S_T=b$ is given by
	\begin{eqnarray*}
		\lefteqn {\p\{M_{H_I,H_{S}}^-<d|H_I<H_S,I_T=a,S_T=b\}}\\
		& \displaystyle{ =
			\frac{  {W}^{(\gamma)}(b-a)}{{W}^{(\gamma)}(d)} e^{-(b-a-d)\left(\frac{{W}^{(\gamma)'}(d)}{{W}^{(\gamma)}(d)} -\Phi(\gamma){W}^{(\gamma)}(d)+\Phi(\gamma)\right)} }
	\end{eqnarray*}
	for $0<d<b-a$.
	\vspace{3mm}
	
	ii. the distribution of the maximum loss of the post-supremum process $M_{H_{S},T}^-$ when $H_I<H_S$, $I_T=a$, $S_T=b$ is given by
	\begin{eqnarray*}\lefteqn{\p(M_{H_{S},T}^-<d|H_I<H_S,I_T=a,S_T=b)} \\
		&  \displaystyle{=\frac{(Z^{(\gamma)}(d)-1)
				\frac{W^{ (\gamma)\prime}(d)}{W^{(\gamma)}(d)} - \gamma W^{(\gamma)}(d) }{  (Z^{(\gamma)}(b-a)-1)
				\frac{W^{ (\gamma)\prime}(b-a)}{W^{(\gamma)}(b-a)} - \gamma W^{(\gamma)}(b-a)   } }
	\end{eqnarray*}
	for $0<d<b-a.$
\end{thm}

\prf
i. For a generic spectrally negative L\'{e}vy process, we have
\begin{eqnarray}
\p_x \{ \tau_{b}^+<\alpha_d, \tau_{b}^+<\tau_{a}^-,\tau_{b}^+<T  \} &=&
\mathbb{E}_{x}\left[e^{-q \tau_{b}^{+}} ; \tau_{b}^{+} \leq   \alpha_{d}\wedge \tau_{a}^{-} \right] \nonumber \\
&=&\frac{{W}^{(\gamma)}(x-a)}{{W}^{(\gamma)}(d)} e^{-(b-d-a)\frac{{{W}}^{(\gamma)'}(d)}{{{W}}^{(\gamma)}(d)}}\label{rev1}
\end{eqnarray}
by the identity \cite[Eq.20]{avram2019w}, top row and middle column, given in terms of the drawdown time $\alpha_d$. In view of \eqref{11} in the proof of Proposition \ref{prop2}, and \eqref{rev1} above, we get
\begin{eqnarray*}
	&& \hspace{-3cm}\p\{M_{H_I,H_{S}}^-<d|H_I<H_S,I_T=a,S_T=b\}\\
	&= & \displaystyle{\lim_{x\rightarrow a}  \frac{\mathbb{P}_{x}^{\Phi(\gamma)}\left\{
			\tau_{b}^+<\alpha_d, \tau_{b}^+<\tau_{a}^-,\tau_{b}^+<T
			\right\}}{ \mathbb{P}_{x}^{\Phi(\gamma)} \{\tau_{b}^+<\tau_{a}^- \wedge T\} } }
	\\	&= &\displaystyle{\lim_{x\rightarrow a} \frac{\T{{W}}^{(\gamma)}(b-a)}{\T{W}^{(\gamma)}(x-a)}\, \frac{\T{W}^{(\gamma)}(x-a)}{\T{W}^{(\gamma)}(d)} e^{-(b-d-a)\frac{\T{{W}}^{(\gamma)'}(d)}{\T{{W}}^{(\gamma)}(d)}} }  \\
	&= &\displaystyle{\frac{\T{W}^{(\gamma)}(b-a)}{\T{W}^{(\gamma)}(d)} e^{-(b-d-a)\frac{\T{W}^{(\gamma)'}(d)}{\T{W}^{(\gamma)}(d)}}  }
	\\ &=& \displaystyle{
		\frac{  {W}^{(\gamma)}(b-a)}{{W}^{(\gamma)}(d)} e^{-(b-d-a)\left(\frac{{W}^{(\gamma)'}(d)}{{W}^{(\gamma)}(d)} -\Phi(\gamma){W}^{(\gamma)}(d)+\Phi(\gamma)\right)} }
\end{eqnarray*}
for $0<d<b-a$, where we used the identity $ \T{W}^{(0)}_{\Phi(\gamma)}(x)= e^{-\Phi(\gamma)x} {W}^{(\gamma)}(x)$.

ii. In view of the conditioning revealed by $h$ in Proposition \ref{prop2}.ii,  we have
\begin{eqnarray*}\lefteqn {\p\{M_{H_{S},T}^-<d|H_I<H_S,I_T=a,S_T=b\}}\\
	&= \displaystyle{ \lim_{x \to b-} \frac{\p_{x}\left\{T<\tau_{b-d}^-, T<{\tau}^+_{b} \wedge{\tau}_a^-  \right\} }{     \p_{x} \{ T<{\tau}^+_{b} \wedge{\tau}_a^-\}          } }\\
	&= \displaystyle{ \lim_{x \to b-} \frac{\p_{x}\left\{T<\tau_{b-d}^-\wedge {\tau}^+_{b} \wedge {\tau}_a^- \right\} }{    h(x)      }    }\\
	&= \displaystyle{ \lim_{x \to b-} \frac{\p_{x}\left\{T<\tau_{b-d}^-\wedge {\tau}^+_{b}  \right\} }{    h(x)      }    }\\
	&= \displaystyle{\lim_{x \to b-}  \frac{1-Z^{(\gamma)}(x-b+d)
			+ (Z^{(\gamma)}(d)-1)\frac{W^{ (\gamma)}(x-b+d)}{W^{(\gamma)}(d)} }{  1-Z^{(\gamma)}(x-a)
			+ (Z^{(\gamma)}(b-a)-1)\frac{W^{ (\gamma)}(x-a)}{W^{(\gamma)}(b-a)}    }    }  \\
	&= \displaystyle{ \frac{(Z^{(\gamma)}(d)-1)
			\frac{W^{ (\gamma)\prime}(d)}{W^{(\gamma)}(d)} -Z^{(\gamma)\prime}(d) }{  (Z^{(\gamma)}(b-a)-1)
			\frac{W^{ (\gamma)\prime}(b-a)}{W^{(\gamma)}(b-a)} -Z^{(\gamma)\prime}(b-a)   }}
\end{eqnarray*}
where we have used $\tau_{b-d}^- < \tau_a^-$ since $a<b-d$, and L'Hospital's rule to take the limit. The result follows by $Z^{(\gamma)\prime}=\Phi(\gamma) W^{(\gamma)}$.
\finish

As a special case, the results of Theorem \ref{thm2} are consistent with the formulas for Brownian motion given in \cite[Prop.4.1]{salm}. The following corollary is immediate for the Laplace transform of the drawdown duration $\mathcal{T}^d_{H_{S},T}$ for the post-supremum process under the conditions of Theorem \ref{thm2}.
\begin{cor} We have
	\begin{eqnarray*}
		\lefteqn{\mathbb{E}[e^{ -\gamma \mathcal{T}^d_{H_{S},T}} | H_I<H_S,I_T=a,S_T=b] }\\
		& \displaystyle{=   1-   \frac{(Z^{(\gamma)}(d)-1)
				\frac{W^{ (\gamma)\prime}(d)}{W^{(\gamma)}(d)} - \gamma W^{(\gamma)}(d) }{  (Z^{(\gamma)}(b-a)-1)
				\frac{W^{ (\gamma)\prime}(b-a)}{W^{(\gamma)}(b-a)} - \gamma W^{(\gamma)}(b-a)   }}
		\; .
	\end{eqnarray*}
	for $0<d<b-a.$
\end{cor}

Similar to the analysis of various distributions found above as $h$-transforms, the distribution of the drawdown duration can be identified as we give next in the following proposition.
In particular in Corollary \ref{cor1}, the law of the drawdown duration of the post-supremum process does not depend on the value of the supremum.  By the same analogy, we now take the process up to $\alpha_d$ rather than an independent exponential time $T$, and consider the part of the process after $\kappa_{\alpha_d}$ until $\alpha_d$. Identification of the underlying conditioning at this part of the path leads to the distribution of the drawdown duration. Indeed the law of the drawdown duration $\mathcal{T}^{d}$ does not depend on the value of $\bar{X}_{\alpha_d}$.

\begin{prop}\label{prop3} Given that $\bar{X}_{\alpha_d}=m$, the law of the post-$\kappa_{\alpha_d}$ process  is a $h$-transform of  the spectrally negative L\'{e}vy process killed at $\alpha_{d}< \tau_m^+ $ with
	\[
	h(x)=1-e^{(m-x)\frac{W^{'}(d)}{W(d)}}\quad  x\le m
	\]
	and
	\[
	\mathbb{E}[e^{-\gamma \mathcal{T}^{d}}]=\frac{W(d)}{W^{'}(d)}\frac{Z^{(\gamma)}(d)W^{(\gamma)'}(d)-\gamma W^{(\gamma)}(d)^2}{W^{(\gamma)}(d)} \; .
	\]
\end{prop}
\prf Given $\bar{X}_{\alpha_{d}}=m,$ the law of the process between $\tau_m^+$  and $\alpha_{d}$ evolves as a spectrally negative L\'evy process conditioned to stay below $m$ and killed at time $\alpha_{d}$. That is, we have
\begin{eqnarray}
& & \hspace{-1cm}\p_{x}\{ {X}_t\in dy, t<\alpha_{d}\wedge {\tau}_m^+ \,|\, \alpha_{d}<   {\tau}_m^+\} \nonumber\\
& =&   \frac{ \p _x  \left\{   X_t\in dy, t<\alpha_{d}  \right\} \p_x\{ \alpha_{d}< \tau_m^+ | X_t = y, t<\alpha_{d}\wedge {\tau}_m^+\}}{\p _x\{\alpha_{d}<\tau_m^+\} }  \nonumber \\
& =& \frac{ \p _x  \left\{   X_t\in dy, t<\alpha_{d} \wedge {\tau}_m^+\ \right\} \p_y\{ \alpha_{d}\circ \theta_t<\tau_m^+\circ \theta_t \}}{\p _x\{\alpha_{d} < \tau_m^+\} } \nonumber \\
& = & \frac{  \p_{y} \{ \alpha_{d}< {\tau}_m^+\} }{\p _{x}\{ \alpha_{d}<  {\tau}_m^+\}} \p_{x}  \{  {X}_t\in dy, t< \alpha_{d} \wedge {\tau}_m^+\}  \nonumber
\end{eqnarray}
by Markov property.
This is the $h$-transform of the law of a spectrally negative L\'{e}vy process killed at $\alpha_{d}$ with
\begin{eqnarray*}
	h(z) &=&   \p_{z} \{\alpha_{d} <\tau_m^+\}  \\
	& =&  1-e^{(m-z)\frac{W^{'}(d)}{W(d)}}
\end{eqnarray*}
by \cite[Thm.1]{arxive}.
Then, the Laplace transform of the drawdown duration can be found using an independent exponential random variable $T$ with rate $\gamma$ for calculation purposes as
\begin{eqnarray}\label{llz}
\mathbb{E}[e^{-\gamma \mathcal{T}^{d}}| \bar{X}_{\alpha_{d}}=m] &=&\lim_{x \to m}P_x\{\alpha_{d}<T|\alpha_{d}<\tau^{+}_{m}\} \nonumber \\
&=& \lim_{x \to  \nonumber m}\frac{P_x\{\alpha_{d}<T,\alpha_{d}<\tau^{+}_{m}\}}{P_x\{\alpha_{d}<\tau^{+}_{m}\}}\\ \nonumber
&=&\lim_{x \to m}\frac{(1-e^{(m-x)\frac{W^{(\gamma)'}(d)}{W^{(\gamma)}(d)}})[Z^{(\gamma)}(d)-W^{(\gamma)}(d)\frac{\gamma W^{(\gamma)}(d)}{W^{(\gamma)'}(d)}]}{1-e^{(m-x)\frac{W^{'}(d)}{W(d)}}}\\\nonumber
&=&  \frac{W(d)}{W^{'}(d)}\frac{Z^{(\gamma)}(d)W^{(\gamma)'}(d)-\gamma W^{(\gamma)}(d)^2}{W^{(\gamma)}(d)}
\end{eqnarray}
where the expression in the numerator is due to  \cite[Eqn.13]{avram2019w} and the denominator follows from \cite[Thm.1]{arxive} (as for the $h$-function) in the last equality. \finish

\begin{Rem} In \cite{LLZ17}, the drawdown duration at time $t\ge 0$ is defined as $t-G_{t}$, where
	$G_{t}:=\sup\{0\leq s\leq t:Y_{s}=0\}$ is the last time the drawdown process is at level $0.$
	The joint distribution of $\alpha_d, G_{\alpha_{d}},Y_{\alpha_{d}},\bar{X}_{\alpha_d}$ is characterized with their Laplace transform for a spectrally negative Levy process considered up to a fixed time $t$ in \cite[Thm 3.1]{LLZ17}. By taking $s=\delta=0,$ $q=\gamma$ and $r=-\gamma,$ the Laplace transform of $\alpha_d-G_{\alpha_{d}}$ is obtained as
	\[
	\mathbb{E} \bigl[e^{-\gamma(\alpha_d-G_{\alpha_{d}})} \bigr]=\frac{W(d)}{W^{'}(d)}\frac{Z^{(\gamma)}(d)W^{(\gamma)'}(d)-\gamma W^{(\gamma)}(d)^2}{W^{(\gamma)}(d)}
	\]
	which matches with the result given in \eqref{llz} as expected  since the drawdown duration $\mathcal{T}^{d}$ is indeed $\alpha_d-G_{\alpha_{d}}$.
\end{Rem}

Finally, for removing the conditions in Theorems \ref{thm3} and \ref{thm2}  through some simple arguments, the joint distribution of the supremum and the infimum can be derived as
\begin{eqnarray*} 
	\!\!  \p_{0}\{a<I_{T},S_{T}<b\} & = & \p_0\{T<\tau^{-}_a\wedge\tau^{+}_b\} \\
	&=& \p_{0}\{T<\tau^{-}_a,\tau^{-}_a<\tau^{+}_b\}+\p_{0}\{T<\tau^{+}_b,\tau^{+}_b<\tau^{-}_a\} \nonumber\\
	&=& 1-(\p_{0}\{T>\tau^{-}_a,\tau^{-}_a<\tau^{+}_b\}+\p_{0}\{T>\tau^{+}_b,\tau^{+}_b<\tau^{-}_a\})\nonumber\\
	&=&1-(\E_{0}[e^{-\gamma\tau^{-}_a}1_{\{\tau^{-}_a<\tau^{+}_b\}}]+\E_{0}[e^{-\gamma\tau^{+}_b}1_{\{\tau^{+}_b<\tau^{-}_a\}}])\nonumber\\
	& =& 1-(\E_{-a}[e^{-\gamma\tau^{-}_0}1_{\{\tau^{-}_0<\tau^{+}_{b-a}\}}]+\E_{-a}[e^{-\gamma\tau^{+}_{b-a}}1_{\{\tau^{+}_{b-a}<\tau^{-}_0\}}]).\nonumber \\
	&=& 1-Z^{(\gamma)}(-a)+[Z^{(\gamma)}(b-a)-1]\frac{W^{(\gamma)}(-a)}{W^{(\gamma)}(b-a)}
\end{eqnarray*}
for $a<0<b$, using the two-sided exit formulas \cite{Kyp}.
This formula for the joint distribution of $I_T$ and $S_T$ of a spectrally negative L\'{e}vy process generalizes that for a Brownian motion with drift \cite{ceren}.

\section{Proof of   Lemma \ref{lem1} \la{s:prf}}
\prf a) Note that $X_{H_{S}-}=S_T$ as $X$ is spectrally negative. Since we have assumed that $X$ is of unbounded variation, which in turn is equivalent to regularity of $0$  for $(-\infty,0)$ \cite[pg.232]{Kyp}. Lem.VI.6.ii of \cite{Ber} for Markov processes yields the independence conditionally on $S_T$.

b) The pre-$H_S$ process evolves as a spectrally negative L\'evy process that stays below $b$ before an exponential time $T$, given that the supremum is reached during the exponential horizon. Its transition function $P_t$ can be identified by
\begin{equation}
P_t(x,dy) = \p_x\{X_t\in dy, t < \tau_b^+\wedge T\,|\, \tau_b^+<T\}, y \leq b  \label{preHS}
\end{equation}
which is equal to
\begin{eqnarray*}
	& &  \displaystyle{\frac{  \p_x\{X_t\in dy, t<T \wedge \tau_b^+, \tau_b^+<T\} }{\p _x\{ \tau_b^+<T\} }} \nonumber \\
	& & =\displaystyle{ \frac{ \p _x  \left\{   X_t\in dy, t<T \wedge \tau_b^+ \right\} \p_x\{ \tau_b^+<T | X_t = y, t<T \wedge \tau_b^+ \}}{\p _x\{\tau_b^+<T\} } } \nonumber \\
	& & = \displaystyle{\frac{ \p _x  \left\{   X_t\in dy, t<T \wedge \tau_b^+ \right\} \p_y\{ T\circ \theta_t>\tau_b^+\circ \theta_t \}}{\p _x\{T> \tau_b^+\} } }\nonumber \\
	& & = \displaystyle{ \frac{  \p_{y} \{ T> \tau_b^+\} }{\p _x\{T> \tau_b^+\}} \p _x  \left\{   X_t\in dy, t<T \wedge \tau_b^+ \right\} }
\end{eqnarray*}
by Markov property and the fact that $T$ is memoryless, where $\theta$ is the shift operator \cite[Sec.7.3]{cinl}. From \eqref{preHS}, we have
\begin{equation}   \label{eqSec6}
P_t(x,dy) = \frac{h(y)}{h(x)}\p _x  \left\{   X_t\in dy, t<T \wedge \tau_b^+ \right\}
\end{equation}
where $h(z)=   \p_{z} \{T>\tau_b^+\}= e^{-\Phi(\gamma)(b-z)}$ by \cite[pg.232]{Kyp}.
Now, note that the entrance law for the pre-supremum process is obtained as $x\to 0$ in \eqref{eqSec6} and is equivalent to the change of measure \eqref{change} because
\begin{eqnarray*}
	\lim_{x\to 0} P_t(x,dy) &=& \lim_{x\to 0} \frac{e^{-\Phi(\gamma)(b-y)}}{e^{-\Phi(\gamma)(b-x)}}\p _x  \left\{   X_t\in dy, t<T \wedge \tau_b^+ \right\} \\
	&=& e^{\Phi(\gamma)y} \, \p   \left\{   X_t\in dy, t<T \wedge \tau_b^+ \right\} \\
	&=& e^{\Phi(\gamma)y-\gamma t} \, \p   \left\{   X_t\in dy, t< \tau_b^+ \right\}
\end{eqnarray*}
which is indeed the Esscher exponential change of measure given in \eqref{change}.
An alternative proof of this result is available  in Theorem 1.ii. of \cite{arxive} based on the excursions of the process from its running maximum. In the expression for the joint distribution of the pre-supremum and post-supremum processes given in \cite[Eq.5]{arxive},    the conditional law of the pre-supremum  appears as a product with the marginal distribution of the supremum, which is  exponential  with parameter $\Phi(\gamma)$.

c) The law of the post-$H_S$ process can be characterized by noting that the post-$H_S$ process evolves as a spectrally negative L\'evy process which is conditioned to stay below $b$ and killed at an exponential time $T$. That is, we have
\begin{eqnarray}
& & \hspace{-1cm}\p_{x}\{ {X}_t\in dy, t<T\wedge \tau_b^+\,|\, T<   {\tau}_b^+\} \nonumber\\
& =&   \frac{ \p _x  \left\{   X_t\in dy, t<T\wedge \tau_b^+  \right\} \p_x\{ T< \tau_b^+ | X_t = y, t<T\wedge \tau_b^+ \}}{\p _x\{T<\tau_b^+\} }  \nonumber \\
& =& \frac{ \p _x  \left\{   X_t\in dy, t<T\wedge \tau_b^+  \right\} \p_y\{ T\circ \theta_t<\tau_b^+\circ \theta_t \}}{\p _x\{T< \tau_b^+\} } \nonumber \\
& = & \frac{  \p_{y} \{T< {\tau}_b^+\} }{\p _{x}\{T<  {\tau}_b^+\}} \p_{x}  \{  {X}_t\in dy, t<T\wedge \tau_b^+  \} \nonumber
\end{eqnarray}
by Markov property and the fact that $T$ is memoryless.
This is the $h$-transform of the law of a spectrally negative L\'{e}vy process killed at the minimum of an exponential time and passage time above $b$, with
\begin{eqnarray}
h(z) &=&   \p_{z} \{T<\tau_b^+\}=\p \{T<\tau_{b-z}^+\} \nonumber \\
& =& 1-\E[e^{-\gamma \tau_{b-z}^+}1_{\{\tau_{b-z}^+<\infty\}}] \nonumber \\
&=& 1-e^{-\Phi(\gamma)(b-z)} \label{ilkh}
\end{eqnarray}
for $z<b$ \cite[pg.232]{Kyp}.
That is, when $S_T=b$, the transition semigroup of the post-$H_S$ process is given by
\[
P_t(x,dy)=\frac{h(y)}{h(x)}\, \p_{x}  \left\{  {X}_t\in dy, t<T \right\}
\]
for $x,y<b$ and the entrance law is obtained as $x\to b-$.
\finish




\vspace{0.5cm}

\noindent{\bf Acknowledgement.}
We are grateful to Andreas E. Kyprianou for
helpful discussions on the post-infimum process.
\vspace{0.5cm}
\\


\begin{thebibliography}{apalike}

\bibitem{arxive} Vardar-Acar C. and Caglar M. 2018. Path Decomposition of Spectrally Negative {L}{\'e}vy Processes,  arXiv preprint arXiv:1801.06426.

\bibitem{mill}  Millar, P. W. 1977.  Zero-one laws and the minimum of a Markov process,Transactions of the American Mathematical Society, 226,365-391.

\bibitem{avram2019w} Avram, F. and Grahovac, D. and Vardar-Acar, C. 2019. The W, Z/$\nu$, $\delta$ Paradigm for the First Passage of Strong Markov Processes without Positive Jumps,  Risks, 7(1), 18-33.

\bibitem{duquesne2003path} Duquesne, T. 2003.
Path decompositions for real L{\'e}vy processes,  Annales de l'IHP Probabilit{\'e}s et statistiques, 39, 2, 339-370.

\bibitem{bertoin1993splitting} Bertoin, J. 1993. Splitting at the infimum and excursions in half-lines for random walks and L{\'e}vy processes,  Stochastic processes and their applications, 47, 1, 17-35.

\bibitem{chaum} Chaumont, L. 1996. Conditionings and path decompositions for L{\'e}vy processes
 Stochastic processes and their applications, 64, 1, 39-54

\bibitem{chaumDon} Chaumont, L. and Doney, R. 2005.
On L{\'e}vy processes conditioned to stay positive,  Electronic Journal of Probability, 10, 948-961

\bibitem{AKP} Avram, F. and Kyprianou, A. and Pistorius, M. 2004. Exit problems for spectrally negative {L}{\'e}vy processes and applications to ({C}anadized) {R}ussian options, The Annals of Applied Probability, 1, 14,  215-238.

\bibitem{ceren} Vardar-Acar, C. and Zirbel, C. L. and Sz{\'e}kely, G. J. 2013.
On the correlation of the supremum and the infimum and of maximum gain and maximum loss of Brownian motion with drift, Journal of Computational and Applied Mathematics, 248, 61-75

\bibitem{AGV}
Avram, F. and Grahovac, D. and Vardar-Acar, C. 2017. 
	The $W,Z$ scale functions kit for first passage problems of  spectrally negative  L\'evy
			processes, and applications to the optimization of dividends,
arXiv preprint arXiv:1706.06841

\bibitem{baur} Baurdoux, E. J. and Palmowski, Z. and Pistorius, M. R. 2017.
On future drawdowns of L{\'e}vy processes, Stochastic Processes and their Applications, 127, 8, 2679-2698.


\bibitem{cinl} \c{C}{\i}nlar, E. 2011. Probability and Stochastics, 261,
Springer Science \& Business Media


\bibitem{KKR} Kuznetsov, A. and Kyprianou, A. E and Rivero, V. 2013. The theory of scale functions for spectrally negative {L}{\'e}vy processes,
{L}{\'e}vy Matters II,
97-186, Springer

\bibitem{Kyp}
Kyprianou, A. 2014. Fluctuations of {L}{\'e}vy Processes with Applications: Introductory Lectures, Springer Science \& Business Media,

\bibitem{LLZ17} Landriault, D. and Li, B. and Zhang, H. 2017. On magnitude, asymptotics and duration of drawdowns for {L}{\'e}vy models,
Bernoulli, 1, 23, 432-458.

\bibitem{MijoPist}
Mijatovic, A. and Pistorius, M. R. 2012. On the drawdown of completely asymmetric {L}{\'e}vy processes, Stochastic Processes and their Applications,
11, 122, 3812-3836.

\bibitem{LLZU} Landriault, D. and Li, B. and Li, S. 2015. Analysis of a drawdown-based regime-switching {L}{\'e}vy insurance model, Insurance: Mathematics and Economics, 60, 98-107.

\bibitem{Ber} Bertoin, J. 1998. {L}{\'e}vy processes, Cambridge university press.

\bibitem{salm} Salminen, P. and Vallois, P. 2007
On maximum increase and decrease of Brownian motion, Annales de l'Institut Henri Poincare (B) Probability and Statistics, 43, 6, 655-676.


	
\end{thebibliography}
\end{document}